# FUZZY CONTINUITY IN SCALABLE TOPOLOGY

**M. Burgin**

Department of Mathematics
*University of California, Los Angeles*
Los Angeles, CA 90095

**Abstract:** In this paper, topological spaces are enriched by additional structures in order to give a more realistic representation of real life phenomena and computational processes and at the same time, to provide for utilization of the powerful technique developed in topology. The suggested approach is based on the concept of a discontinuity structure or scale Q in a topological space *X*. This structure serves as a scale for the initial topology on *X*. Problems of science and engineering need such a scale because all measurements and the majority of computations are performed in a definite scale and problems of scalability are important both for physical theories and computational algorithms. Taking a mapping of a topological space *X* with the discontinuity structure Q into a topological space *Y* with a discontinuity structure R, we define (Q,R)-continuity, weak (Q,R)-continuity, and R-continuity of this mapping. Fuzzy continuous functions, which are studied in neoclassical analysis, are examples of R-continuous mappings. Different properties of scales, Q-open, Q-closed sets, and (Q,R)-continuous mappings are obtained.

# 1. Introduction

Topology is becoming one of the main mathematical tools in physics. Even when physicists work with real and complex numbers, they utilize topological properties of these numbers. To find derivatives and integrals, one needs existence of limits. In turn, limits are topological constructions. It means that manipulations with limits and, consequently, with derivatives and integrals are based on the topology in number spaces. Number spaces, the real line and complex space, as well as Euclidean spaces have a good topological structure that allows mathematicians to develop calculus and optimization methods in these spaces. Such topological structures are called metric spaces and have a lot of useful features. These features provide for solution of many theoretical and practical problems. All this shows that topology is an important field of mathematics and important tool for science, especially, for contemporary physics where rather abstract topological spaces have become efficient means for modeling and researching (cf., for example, (Witten, 1987; 1988; Nash, 1997)). Such fields as topological quantum field theory become more and more influential.

Now topological technique also finds its applications in computer science (cf., for example, (Burgin, 2005)). Topology comes to computations from different directions.

First, numerical methods are built on topological and metric properties of numbers. Results are produced with some precision, which estimates the distance to the completely precise result, which is possible, in general, only as a theoretical abstraction.

Second, computers are used for working with topological models. Simulation of different physical and technological processes gives many examples of such situations as we know that topology is inherent to physics as physical models use more and more topology. Simulation of industrial processes and advanced measuring techniques allow experts to penetrate deeply into the details of flow and transport phenomena as occurring in industrial equipment, and to analyze their mutual relations, properties, and their impact on physical and chemical processes. Simulation of sub-micron devices and multi-length scale materials phenomena allows physicists to get new knowledge and engineers to build new powerful devices.

Third, models of computations themselves involve topological structures. For instance computational trajectories form a topological space (Burgin, 2001a). In many cases, properties of this space determine results of computation (cf., for example, (Büchi, 1960; Muller, 1963; Burgin, 1992; Vardi and Volper, 1994; Boldi and Vigna, 1998)). For instance, a finite automaton *A* accepts an infinite string when *A* comes to some accepting state is an adherent point of the computational trajectory. Such automata are called Büchi automata (Büchi, 1960). Now there is a developed theory of finite automata working with infinite words (cf., for example, (Vardi and Volper, 1994). The result of a limit Turing machine is the limit of partial results obtained in a computation (Burgin, 1992). One more topological model of computation introduced in (Chadzelek and Hotz, 1999) is δ–***Q***-*machine* which takes in real numbers as inputs and use infinite converging computations on more and more precise rational roundings of their inputs.

One more peculiarity of contemporary science manifests itself in multi-scale physics. This phenomenon is rooted both in real life situations and in intrinsic structures of our universe. Let us consider at first real life situations.

When we look at some polished body, such as a table or mirror, we do not see holes and the surface looks smooth and continuous. However, taking a microscope, we notice a lot of holes and essential unevenness. Thus, what is smooth and continuous on the scale of our natural vision becomes uneven and discontinuous on the microscopic scale. Going deeper to the level of molecules and atoms, we find that any substance consists mostly of holes.

At the same time, when we closely look at an asphalt road, we see that it is uneven and has many holes. However, when we walk, we do not notice these holes. In other words, the asphalt surface is smooth and continuous for walking, but does not look so from a close distance.

Moreover, when walking, we can ignore holes with the diameter 0.1 in, but cannot ignore holes with the diameter 10 in. We can overstep holes with the diameter 10 in, but have to go around holes with the diameter 10 ft.

All these examples show that in contrast to mathematics, in real life what is considered smooth and continuous depends on the scale on which we operate.

Moreover, in a similar way, scalability emerges in various fields of physics. From the historical perspective, the most apparent case of scale dependency was discovered at the beginning of the 20$^{th}$ century. It was demonstrated that classical physics that for centuries worked well for planets and bodies on the Earth was not adequate for very small length scales, as well as for very large velocity scales. In the first case, quantum mechanics substitutes Newtonian mechanics as a more adequate theory, while in the second case, relativity theory becomes necessary to achieve sufficient theoretical precision with respect to experimental data.

Now when physics becomes more and more sophisticated, there is an active discussion whether our space-time is continuous or discrete. As Motl writes, "the idea that space-time could be discrete has been a recurring one in the scientific discourse of the twentieth century. A survey of just a few examples reveals that discrete space-time can actually mean many things and is motivated by a variety of philosophical or theoretical influences.

It has been apparent since early times that there is something different between the mathematical properties of the real numbers and the quantities of measurement in physics at small scales. Riemann himself remarked on this disparity even as he constructed the formalism, which would be used to describe the space-time continuum for the next century of physics. When you measure a distance or time interval you can not declare the result to be rational or irrational no matter how accurate you manage to be. Furthermore, it appears that there is a limit to the amount of detail contained in a volume of space."

Physicists show that if they try to use an interferometer, or simple time of flight measurements to determine locality, they get the answer that the minimal distance measureable is the Planck length. So, there really is a sense in which distance shorter than the Planck length has no meaning (Calmet, Graesser, and Hsu, 2004). This supports the hypothesis in physics there is a minimal length interval and minimal time interval connected to the Planck constant that space, time, and matter all have to be discrete ('t Hooft, 1996; 1999).

These problems brought science to Planck scale physics (cf., for example, (Ring, *et al*, 1995; Requardt, 1998; Albrecht and Skordis, 2000; Brandenberger and Martin, 2001) or as some call it, string scale physics. It studies physical phenomena in space-time domains determined by the fundamental length and time intervals. In the domains at this scale, named after Max Planck, one of the founders of quantum theory, many properties of physical space and time are essentially different. As Requardt writes (1998), "starting from the working hypothesis that both physics and the corresponding mathematics have to be described by means of discrete concepts on the Planck scale, one of the many problems one has to face in this enterprise is to find the discrete protoforms."

Physicists even go deeper to the super-Planck scale physics. For instance, Brandenberger and Martin (2001) calculate the spectrum of density fluctuations in models of inflation based on a weakly self-coupled scalar matter field minimally coupled to gravity, and specifically investigate the dependence of the predictions on modifications of the physics on length scales smaller than the Planck length.

At the same time, the development of nanotechnology in general and nanoelectronics, in particular, is based on nanoscale physics (cf., for example, (Parfitt, 1997; Kirsch, 2004; Riley, 2004)). The flow of research in this area is so active that a new journal " Virtual Journal of Nanoscale Science & Technology" has been organized.

In addition, physicists study small-scale and high energy scale phenomena (cf., for example, (Batchelor, 1959; Belyayev, et al, 1974; Caldwell, et al, 1975; Caldwell, 1983; Ring, *et al*, 1995)), fine scale structures (cf., for example, (Caldwell, 1976), develop sub-grid physics (Porté-Agel, et al, 2001; 2001a), reactor-scale burning plasma physics and meso-scale physics, and apply large scale physics (cf., for example, (Beazley and Lomdahl, 1997)). Thus, physics starts taking into account a variety of scales even in one area. For instance, many phenomena in materials involve atomic motion on a wide range of length scales from the atomic to the mesoscopic to the macroscopic. Turbulence and crack propagation are well-known examples. Change of the scale demands transformation or even change of the used model. For instance, as the length scale increases only a decreasing fraction of the information contained in the motion of each atom is relevant to many phenomena, and it becomes unnecessary to describe the motion of all atoms on the atomic scale. As a result, we come to mutli-scale physics.

A similar or may be, even more urgent situation with scalability exists in the field of information processing. Scalability is a highly significant issue in computer simulation, electronics systems, data and knowledge bases, networking, etc.

At the same time in classical mathematics, topology in general and its principal concept continuity are not scalable. For example, a mapping is either continuous or not and its continuity does not depend on any scale. To eliminate this restriction and to introduce scalability into topology, we need additional structures in topological spaces with different properties to reflect more adequately computational and physical reality than it is done by classical structures in analysis and topology. Such new structures together with relevant methods and constructions are provided by *scalable* or *discontinuous topology* (Burgin, 2004). It is a new field in which ordinary structures of topology are studied by means of scales and in some sense, fuzzy concepts. For example, the continuous mappings studied in classical topology become a part of the set of the *scaled continuous or fuzzy continuous mappings* studied in scalable topology (Burgin, 2004). The scale of a scaled topological space defines to what extent it is acceptable to extend continuity.

According to Lefschetz (1930), "Topology or Analysis Situs is usually defined as the study of properties of spaces or their configurations invariant under continuous transformations. But, Lefschetz continues, what are spaces and continuous transformations?" This is a crucial question for topology and answer to this question defines the essence of topology. That is why one of the topologists writes that topology emerged and is growing in a series of experiments. A new step in this direction is introduction of scalable topology. Its main idea is to take conventional topological spaces as initial spaces but to relax conditions on continuous mappings. Scales show what discontinuities are tolerated. It makes possible to enlarge the family of acceptable transformations, allowing, for example, transformations in which breaks of continuity are not very large or not very frequent or we do not have precise knowledge about these breaks.

This new version of topology is called discontinuous because the main emphasis is on discontinuous mappings, while topological categories are built with respect to discontinuous mappings. It is also called scalable because discontinuity of mappings and architecture of other topological structures in this theory are determined by a scale defined in a topological space. This approach is different from the classical topology where only continuous mappings are considered. At the same time, the methods that are used for study of these discontinuous mappings are developed from the methods of the classical topology. Such an approach extends the scope of topology making, at the same time, its methods more precise in many situations, especially, in applications. Consequently, new results are obtained extending and even completing classical theorems.

There were several attempts to introduce mathematical structures that are similar but less restrictive than the conventional topology (Appert, and Ky Fan, 1951; Cech, 1968; Hammer, 1962; Netzer, 1978; Sierpinski, 1934). The aim of discontinuous topology is different. A new structure called a *discontinuous structure* or *scale* is added to a conventional topological space as an approximation to the space topology. Its properties of a scale manifest in permissible deviance from standard topological structures – a less precise scale allows more deviations than a more precise one. For instance, the scale determines what size/kind of discontinuity we can ignore. This new structure defines a *distructured* or *scaled topological space* $(X, T_X, Q_X)$ where the pair $(X, T_X)$ is a conventional topological space $X$ with a topology $T_X$. The scale $Q_X$ is some weakening of $T_X$. It means that the initial topology $T_X$ on $X$ is not ignored. It remains as a reference frame for the scale $Q_X$ as well as a limit for a family of such scales $Q_X$ on $X$. At the same time, topological spaces may be considered as particular cases of scaled topological spaces in which the scale coincides with the initial topology. The construction of a scaled topological space reflects situations in real life when exact mathematical models are transformed by approximate results of measurement and computations. In such a way, the topology $T_X$ and the scale $Q_X$ complement each other. Namely, the conventional topology $T_X$ being an abstract construction has a highly developed mathematical apparatus, while the scale $Q_X$, being less conventional but more realistic, affords means for a more accurate representation of reality, linking it to a traditional topological system. Topological spaces become scalable, and choice of a scale regulates properties of these spaces and mappings between them.

The development of scalable topology has been initiated in the framework of neoclassical analysis (Burgin, 1995; 2000; 2001). It is a new direction in mathematical analysis. In it, ordinary structures of analysis, that is, functions and operators, are studied by means of fuzzy concepts. For example, the concept of a limit studied in classical analysis becomes a specification of the concept of a fuzzy limit studied in neoclassical analysis. It extends the scope of analysis making, at the same time, its methods more precise. Consequently, new results are obtained extending and even completing many classical theorems. In addition to this, facilities of analytical methods for various applications become more broad and efficient. It is necessary to remark that not all properties of classical constructions, such as limits or continuous functions, remain true for their neoclassical extensions.

In this work scaled topological spaces and their fuzzy continuous mappings are introduced and studied. There are many such mappings in different fields of mathematics. As an example, we can take step functions, which are basic in the theory of integration (Saks, 1937) or membership functions, which are even more important for set theory (Bourbaki, 1960; Fraenkel and Bar-Hillel, 1958). Providing a possibility to investigate such functions by topological methods helps to achieve better understanding of mathematical structures.

At the same time, almost all functions outside mathematics are fuzzy continuous. For example, all computable functions are fuzzy continuous because it is possible to carry out computations only to a definite precision and we have to bear in mind the results of truncation. In other words, any numerical method that does not take into account the roundoff effect may be very misleading by giving an insufficient approximation or even diverging when computer realizes it (Alefeld and Grigorieff, 1980; Burgin and Westman, 2000). Truncation and roundoff operations translate real numbers into rational numbers and replace continuous functions of mathematical models by fuzzy continuous functions that are processed by computers. Thus, we come to a necessity to study fuzzy continuous functions and operations performed with these functions.

It is necessary to remark that those definitions that are given in this paper are not purely formal constructions but introduce mathematical models for different important real phenomena and systems.

In Conclusion, some directions in scalable topology are formulated for future research.

**Denotations:**

***R*** is the real line;

***R***$^+$ is the set of all non-negative real numbers;

***R***$^{++}$ is the set of all positive real numbers;

***R***$^n$ is an *n*-dimensional Euclidean space;

**N** is the set of all natural numbers;

ω is the sequence of all natural numbers;

∅ is the empty set;

If *X* is a set, then $2^X$ denotes the set of all subsets of *X*;

If *X* is a topological space, then T$_X$ denotes the topology in *X*;

O*x* denotes a neighborhood of a point *x*;

If *a* is a positive number and *X* is a metric space, then O$_a$*x* = { $z \in X$; $\rho(z, x) < a$ };

If *C* is a set in a topological space in *X*, then $\overline{C}$ denotes the closure of *C* in *X*;

If *X* and *Y* are topological spaces, then *F*(*X*, *Y*) is the set of all and *C*(*X*, *Y*) is the set of all continuous mappings from *X* into *Y*;

A projection *f*: *X*→ *Y* is a mapping of *X* onto *Y*;

As mappings are special kinds of binary relations (cf., for example, (Bourbaki, 1960)), it is possible to define for them all set theoretical relations like inclusion or intersection without changing denotations.

**2. Scaled Spaces and Fuzzy Continuous Mappings**

We start with conventional topological spaces. Let *X* be a topological space with a topology T$_X$ and TQ be a subset of T$_X$, i.e., TQ consists of open sets from *X*.

**Definition 1.** A *scale* or *discontinuity structure* $Q_X = Q$ on $X$ is a mapping $X \to 2^{TQ}$ that satisfies the following conditions:

**(SC1)** For all points $x$ from $X$, if $A$ is an element of $Q(x)$, then $x \in A$.

**(SC2)** Any set from $TQ$ belongs to some set from $Q(x)$.

**Remark 1.** In some cases, it is natural to assume that scales $Q_X = Q$ on $X$ satisfy an additional condition (**F**): $X$ is an element of $Q(x)$ for all $x$ from $X$. This condition implies axiom (SC2).

Q is called a discontinuity structure because it determines admissible discontinuity, that is, it determines to what extent a mapping from one topological space into another may be discontinuous. Informally, $Q_X$ relates each point $x$ of $X$ to some set of neighborhoods of $x$. For simplicity, neighborhoods from the set $Q(x)$, we call Q-*neighborhoods* of $x$ and denote by Q the set $Q(X)$.

Mappings are special kinds of binary relations (cf., for example, (Bourbaki, 1960)). Consequently, it is possible to define for them all set theoretical relations like inclusion or intersection without changing denotations. If we take set-valued mapping defined on one set $X$, then inclusion of such mappings means inclusion of images for all points of $X$. For example, if Q and P are discontinuity structures on a topological space $X$, then $P \subseteq Q$ means that $P(x) \subseteq Q(x)$ for all $x$ from $X$.

**Remark 2.** It is possible to consider a scale $Q_X$ as nonhomogenious fiber bundle with the base $X$ and fibers $Q(x)$ (Goldblatt, 1979).

**Definition 2.** A *scaled* or *distructured topological space* is a triad $(X, T_X, Q_X)$ where $T_X$ is a topology and $Q_X$ is a scale/discontinuity structure on $X$.

We know that it is possible to define different topologies on an infinite set. In a similar way, it is possible to define different scales in one topological space.

**Example 1.** Let $X$ be an *n*-dimensional Euclidean space $\mathbf{R}^n$ and $a \in \mathbf{R}^{++}$. Then each collection of sets $Q_X(x)$ consists of all neighborhoods of a point $x$ from $X$ that contain a closed *n*-dimensional ball with the radius $a$ and center $x$. We denote such scale $Q_X$ by $Q_a$.

**Example 2.** Let $X$ be an $n$-dimensional Euclidean space $\boldsymbol{R}^n$ and $a \in \boldsymbol{R}^{++}$. Then each collection of sets $Q_X(x)$ consists of all open $n$-dimensional balls with the radius $r > a$ and center $x$. We denote such scale $Q_X$ by $CQ_a$.

**Example 3.** Let $X$ be an $n$-dimensional Euclidean space $\boldsymbol{R}^n$ and $a \in \boldsymbol{R}^{++}$. Then each collection of sets $Q_X(x)$ consists of all neighborhoods of a point $x$ from $X$ that contain an open $n$-dimensional ball with the radius $a$ and center $x$. We denote such scale $Q_X$ by $Q_{Oa}$.

**Example 4.** Let $X$ be an $n$-dimensional Euclidean space $\boldsymbol{R}^n$ and $a \in \boldsymbol{R}^{++}$. Then each collection of sets $Q_X(x)$ consists of all open $n$-dimensional balls with the radius $r \geq a$ and center $x$. We denote such scale $Q_X$ by $CQ_{Oa}$.

**Example 5.** Let $X$ be a metric space with the metric $\rho$ and $a \in \boldsymbol{R}^+$. Then $Q_X(x)$ consists of the space $X$ and all connected neighborhoods of a point $x$ from $X$ containing a closed ball with the radius $a$ and center $x$. We denote such scale $Q_X$ by $Q_a$.

**Example 6.** Let $X$ be a metric space with the metric $\rho$ and $a \in \boldsymbol{R}^+$. Then $Q_X(x)$ consists of all connected neighborhoods of a point $x$ from $X$ containing an open ball with the radius $a$ and center $x$. We denote such scale $Q_X$ by $CQ_a$.

**Example 7 (P-structure).** Let us fix for each point $x$ from $\boldsymbol{R}$ some neighborhood $Ox$. Then we define $Q_X(x)$ as the set of all open sets in $\boldsymbol{R}$ containing $Ox$. Such scale is called a P-structure. The scale $Q_{Oa}$ is a P-structure.

**Definition 3 (F-structure).** For all $x$, the set $Q(x)$ is a filter in $T_X$ of neighborhoods of $x$, that is, $Q(x)$ is closed with respect to finite intersections and supersets of its elements. Scales from the examples 1, 3 are F-structures. Moreover, they are principal F-structures or P-structures because all sets $Q(x)$ are principal filters.

**Remark 3.** Any discontinuity P-structure is a discontinuity F-structure.

**Definition 4.** a) **U-structure:** For all $x$, $Q(x)$ is closed with respect to arbitrary unions with elements from TQ.

b) **weak U-structure:** TQ is closed with respect to arbitrary unions.

**Remark 4.** Any discontinuity F-structure is a discontinuity U-structure.

**Example 8.** Let $X$ be a metric space with the metric $\rho$ and $a \in \mathbf{R}^+$. Then $Q_X(x)$ consists of $X$ and all bounded neighborhoods of a point $x$ from $X$ containing an open ball with the radius $a$ and center $x$. We denote such scale $Q_X$ by $BQ_{Oa}$. It is a U-structure but not an F-structure.

**Definition 5.** a) **I-structure:** For all $x$, $Q(x)$ is closed with respect to the intersection with an arbitrary element from TQ that contains $x$.

b) **weak I-structure:** TQ is closed with respect to arbitrary finite intersections.

**Definition 6.** a) **L-structure:** For all $x$, $Q(x)$ is a lattice of neighborhoods of $x$.

b) **weak L-structure:** TQ is a lattice.

**Example 9.** $Q_X = T_X$, i.e., $Q_X$ relates each point $x$ of $X$ to the set of all neighborhoods of $x$. In this case, $T_X$ is called the trivial scale.

**Definition 7.** A subset H of a scale Q is called a subscale of the scale Q.

**Example 10.** Let $X$ be an $n$-dimensional Euclidean space $\mathbf{R}^n$ and $a, b \in \mathbf{R}^{++}$. The scale $Q_a$ is a subscale of the scale $Q_{Oa}$. If $a < b$, then the scale $Q_{Ob}$ is a subscale of the scale $Q_{Oa}$ and the scale $Q_b$ is a subscale of the scale $Q_a$.

**Definition 8.** Sets from the scale $Q_X$ are called $Q_X$–*open*, while their complements in $X$ are called $Q_X$–*closed*.

Some properties of scaled topological constructions are similar to the corresponding properties of conventional or absolute topological constructions, while other properties are essentially different. For instance, there are no bounded $BQ_{Oa}$-closed sets but the empty set $\varnothing$ with respect to the discontinuity structure $BQ_{Oa}$ in an unbounded metric space $X$ (cf., Example 2.6), e.g., $\mathbf{R}^n$. The empty set $\varnothing$ is not $BQ_{Oa}$-open in $X$. Likewise in a general case, it is possible that the union of two or more $Q_X$–open sets is not $Q_X$–open or the intersection of two or more $Q_X$–open sets is not $Q_X$–open. This is different from properties of conventional topological spaces.

**Proposition 1.** a) A scale $Q_X$ is a weak U-structure if and only if the intersection of Q-closed sets is a Q-closed set.

b) A scale $Q_X$ is a weak I-structure if and only if any finite union of Q-closed sets is a Q-closed set.

**Corollary 1.** A scale $Q_X$ is a weak L-structure if and only if the intersection of two Q-closed sets is a Q-closed set and the union of two Q-closed sets is a Q-closed set.

It is possible to consider operations on scales: binary operations such as union, intersection, etc.; unary operations such as U-closure, I-closure, L-closure, F-closure, etc.

Let us assume that $(X, T_X, Q_X)$ and $(Y, T_Y, R_Y)$ are scaled topological spaces. In what follows, we omit subscripts $X$ and $Y$ from $T_X$, $T_Y$, $Q_X$ and $R_Y$, when it does not lead to confusion.

**Definition 9.** A mapping $f: X \to Y$ is called:
a) (Q, R)-*continuous at a point* $x$ from $X$ if for $y = f(x)$ and any neighborhood $Oy$ from $R(y)$, the set $f^{-1}(Oy)$ is an element of $Q(x)$.
b) *locally* (Q, R)-*continuous* if it is (Q, R)-continuous at all points from $X$.
c) (Q, R)-*continuous* if for any R-open set $V$, the set $f^{-1}(V)$ is Q-open.
d) R-*continuous at a point* $x$ from $X$ if it is $(T_X, R)$-continuous at $x$.
e) R-*continuous* if it is $(T_X, R)$-continuous.
f) *locally* R-*continuous* if it is locally $(T_X, R)$-continuous.

$(Q_X, R_Y)$-continuous mappings may be considered as mappings that are continuous only to some extent. That is why they are called *fuzzy* or *scaled continuous*.

**Example 11.** When $Y$ is a metric space and R is defined as in Example 2.2, i.e., $R(x) = Q_a(x)$ for some $a \in \mathbf{R}^{++}$ and all $x$ from $X$, then all R-continuous mappings from $X$ to $Y$ coincide with the $a$-fuzzy continuous mappings, which are defined in (Burgin, 1997).

**Lemma 1.** A mapping $f: X \to Y$ is continuous (at a point $x$ from $X$) if and only if $f$ is $(T_X, T_Y)$-continuous (at a point $x$ from $X$).

This result shows that the concept of the fuzzy or scaled continuity is a natural extension of the concept of the conventional continuity. In some sense (cf., for example, Theorem 5), the conventional continuity is a limit case of scaled continuities.

**Lemma 2.** A mapping $f: X \to Y$ is continuous if and only if $f$ is locally $(T_X, T_Y)$-continuous.

This result shows that the concept of the local fuzzy or scaled continuity is also a natural extension of the concept of the conventional continuity. The reason is that local and global continuities coincide in the conventional case. Thus, we come to a question whether the same is true for scaled continuities in a general case.

**Proposition 2.** If $f: X \to Y$ is a locally (Q, R)-continuous projection, then $f$ is (Q, R)-continuous.

Indeed, let $f$ be a locally (Q, R)-continuous projection and $V$ be an R-open set. Then $V$ is an R-neighborhood of some point $y$ from $Y$. As $f$ is a projection, there is a point $x$ from $X$ such that $f(x) = y$. By Definitions 9 a) and b), $f^{-1}(V)$ is an element of $Q(x)$, i.e., $f^{-1}(Oy)$ is an Q-open set. As $V$ is an arbitrary R-open set, this concludes the proof of Proposition 2.

It is possible to ask the question whether this result is true for arbitrary mapping. The following example shows that this is not the case.

**Example 12.** Let us take $X = [0, \tfrac{1}{2}) \cup (\tfrac{1}{2}, 1]$, $Y = [0, 1]$ and the mapping $f: X \to Y$ that is identical on $X$, i.e., $f(x) = x$ for all $x \in X$. For an arbitrary point $z$ either from $X$ or from $Y$, both scales $Q(z)$ and $R(z)$ consist of all open symmetric intervals with the center in $z$, i.e., neighborhoods of $z$ have the form $(z - r, z + r)$ where $0 \leq z - r < z + r \leq 1$. This mapping $f$ is locally (Q, R)-continuous, but it is not (Q, R)-continuous because $(0, 1)$ is an R-open set in $Y$, but $f^{-1}(0, 1) = (0, \tfrac{1}{2}) \cup (\tfrac{1}{2}, 1)$ is not a Q-open set in $X$.

It is possible to ask the question whether any (Q, R)-continuous mapping $f: X \to Y$ is locally (Q, R)-continuous. The answer is negative as the following example demonstrates.

**Example 13.** Let us take $X = Y = \mathbf{R}$ with the natural topology. For any rational point $x$ from $X$, $Q_X(x)$ consists of all intervals with rational ends that contains $x$. For any irrational point $x$ from $X$, $Q_X(x)$ consists of all intervals with irrational ends that contains $x$. For any rational point $y$ from $Y$, $R_Y(y)$ consists of all intervals with irrational ends that contains $x$. For any irrational point $y$ from $Y$, $R_Y(y)$ consists of all intervals with rational ends that contains $x$. Both Q and R are discontinuity C-structures in $\mathbf{R}$.

By Definition 9, the identity mapping $\mathbf{1}: X \to Y$ is (Q, R)-continuous (globally) but it is not (Q, R)-continuous at any point of $X$.

Thus, local and global cases do not coincide for scaled continuity. However, in some cases, local and global (Q, R)-continuity do coincide.

Properties of the system of all open sets in a topological space (Kuratowski, 1966) show that the following result is valid.

**Proposition 3.** A mapping $f: X \to Y$ is R-continuous if and only if it is locally R-continuous.

As in the case of the conventional continuity, it is possible to characterize the scaled continuity by conditions on closed sets.

**Proposition 4.** A mapping $f: X \to Y$ is (Q, R)-continuous if and only if the inverse image $f^{-1}(Z)$ of any R-closed in $Y$ set $Z$ is Q-closed in $X$.

This implies the following well-known result (cf., for example, (Kelly, 1957) or (Kuratowski, 1966)).

**Corollary 2.** A mapping $f: X \to Y$ is continuous if and only if the inverse image $f^{-1}(Z)$ of any closed in $Y$ set $Z$ is closed in $X$.

If operations of the intersection and union are defined in the collection ***L*** of all discontinuity structures on a topological space $X$ with a fixed topology, then it transforms ***L*** into a lattice. Consequently, it provides for the construction of ***L***-fuzzy set of continuous mappings. In the case of metric spaces, it is possible to map ***L*** onto the interval [0,1]. Thus, we obtain a fuzzy set of continuous mappings, which is studied in (Burgin, 1995).

Let us consider four scaled topological spaces $(X, T_X, Q)$, $(Z, T_Z, P)$, $(Y, T_Y, H)$, and $(Y, T_Y, R)$ and two mappings $f: X \to Y$ and $g: Y \to Z$. In addition, we assume that $TR \subseteq TH$, i.e., any R-neighborhood is an H-neighborhood.

**Theorem 1.** If the mapping $f$ is (Q, H)-continuous at a point $x$ (locally (Q, H)-continuous, (Q, H)-continuous) and the mapping $g$ is (R, P)-continuous at the point $f(x)$ (locally (R, P)-continuous, (R, P)-continuous), then the mapping $gf$ is (Q, P)-continuous at a point $x$ (locally (Q, P)-continuous, (Q, P)-continuous, correspondingly).

Indeed, let $z \in Z$, $gf(x) = z$, and $Oz$ is a P-neighborhood of $z$. Then $g^{-1}(Oz)$ belongs to R as the mapping $g$ is (R, P)-continuous at the point $f(x)$. Consequently, $g^{-1}(Oz)$ belongs to H as TR $\subseteq$ TH. Thus, $f^{-1}(g^{-1}(Oz))$ belongs to Q as the mapping $f$ is (Q, H)-continuous at the point $x$. It means that the mapping $gf$ is (Q, P)-continuous at a point $x$ as $f^{-1}(g^{-1}(Oz)) = (gf)^{-1}(Oz))$.

In a similar way, we establish local (Q, P)-continuity and (Q, P)-continuity of the mapping $gf$.

**Theorem 2.** If the mapping $f$ is (locally) (Q, H)-continuous continuous, and the mapping $g$ is (locally) (R, P)-continuous at the point $f(x)$, then the mapping $gf$ is (locally) (Q, P)-continuous.

**Corollary 3.** If the mapping $f$ is (Q, R)-continuous at a point $x$ (locally (Q, R)-continuous, (Q, R)-continuous) and the mapping $g$ is (R, P)-continuous at the point $f(x)$ (locally (R, P)-continuous, (R, P)-continuous), then the mapping $gf$ is (Q, P)-continuous at a point $x$ (locally (Q, P)-continuous, (Q, P)-continuous, correspondingly).

Corollaries 2 and 3 allows one to define the following categories:
- The category **SCTOP** of scaled topological spaces and their scaled continuous mappings.
- The category l**SCTOP** of scaled topological spaces and their locally scaled continuous mappings.
- The category **SCTOP**$_x$ of scaled topological spaces and their mappings scaled continuous at a chosen point such the chosen point from the domain is mapped on the chosen point of the range.

**Corollary 4.** If the mapping $f$ is R-continuous at a point $x$ (R-continuous) and the mapping $g$ is (R, P)-continuous at the point $f(x)$ (locally (R, P)-continuous, (R, P)-continuous), then the mapping $gf$ is P-continuous at a point $x$ (P-continuous).

**Corollary 5.** If the mapping $f$ is continuous at a point $x$ (continuous) and the mapping $g$ is P-continuous at the point $f(x)$ (P-continuous), then the mapping $gf$ is P-continuous at a point $x$ (P-continuous).

This implies the following well-known basic topological result.

**Corollary 6.** If the mapping $f$ is continuous at a point $x$ (continuous) and the mapping $g$ is continuous at the point $f(x)$ (continuous), then the mapping $gf$ is P-continuous at a point $x$ (continuous).

**Definition 10.** A mapping $f: X \to Y$ is called:

a) *weakly (Q, R)-continuous* at a point $x$ from $X$ if for $y = f(x)$ and any neighborhood $Oy$ from $R(y)$, there is a neighborhood $Ox$ from $Q(x)$ such that $f(Ox)$ is a subset of $Oy$.

b) *locally weakly (Q, R)-continuous* if it is (Q, R)-continuous at all points from $X$.

c) *weakly (Q, R)-continuous* if for any R-open set $U$, there is a Q-open set $V$ such that $f(V)$ is a subset of $U$.

d) *weakly R-continuous at a point $x$* from $X$ if it is weakly $(T_X, R)$-continuous at $x$.

e) *weakly R-continuous* if it is weakly $(T_X, R)$-continuous.

f) *locally weakly R-continuous* if it is locally weakly $(T_X, R)$-continuous.

Weakly $(Q_X, R_Y)$-continuous mappings may be considered as mappings that are continuous only to some extent. That is why they are called *weakly fuzzy* or *scaled continuous*.

**Example 14.** When $Y$ is a metric space and the scale R in $Y$ is defined as in Example 2, i.e., $R = Q_a$ for some $a$ from $X$, then all weakly $(T_X, Q_a)$-continuous mappings from $X$ to $Y$ coincide with the $a$-fuzzy continuous mappings, which are defined in (Burgin, 1997).

**Lemma 3.** Any (Q, R)-continuous (at a point $x$ from $X$) is weakly (Q, R)-continuous (at $x$).

**Corollary 7.** Any R-continuous (at a point $x$ from $X$) is locally weakly R-continuous (at $x$).

**Corollary 8.** Any locally (Q, R)-continuous is locally weakly (Q, R)-continuous.

**Corollary 9.** Any locally R-continuous is locally weakly R-continuous.

The concept of weak R-continuity and R-continuity are very close as we have the following result.

Let us consider a topological space $(X, T_X)$ and a scaled topological space $(Y, T_Y, R)$ in which the scale R is closed with respect to neighborhoods in the following sense: if $z \in U$ and $U$ is R-open, then there is an R-neighborhood $V$ of the point $y$ such that $V \subseteq U$.

**Lemma 4.** Any (locally) weakly $(T_X, R)$-continuous is (locally) $(T_X, R)$-continuous.

Proof. Let us consider a locally weakly ($T_X$, R)-continuous mapping $f$: $X \to Y$, a point $x \in X$ with $f(x) = y$, and an R-neighborhood $V$ of the point $y$. To prove Lemma 2.4, we need to show that $f^{-1}(V)$ is an open set in $X$.

Let us take a point $v$ from $f^{-1}(V)$ and its image $z = f(v)$. By the definition of a scale, $V$ is an open set in $Y$. Consequently, there is an R-neighborhood $Oz$ of the point $z$ such that $Oz \subseteq V$. By Definition 10, there is a neighborhood $Ov$ of the point $v$ such that $f(Ov) \subseteq Oz$. Thus, $Ov \subseteq f^{-1}(V)$. As $v$ is an arbitrary point from $f^{-1}(V)$, the set $f^{-1}(V)$ is open in $X$. Local continuity is proved.

The proof for global continuity is similar.

**Lemma 5.** A mapping $f$: $X \to Y$ is continuous (at a point $x$ from $X$) if and only if $f$ is weakly ($T_X$, $T_Y$)-continuous (at a point $x$ from $X$).

This result shows that the concept of the weak fuzzy or scaled continuity is also a natural extension of the concept of the conventional continuity.

**Lemma 6.** A mapping $f$: $X \to Y$ is continuous if and only if $f$ is locally weakly ($T_X$, $T_Y$)-continuous.

This result shows that the concept of the local weak fuzzy or scaled continuity is also a natural extension of the concept of the conventional continuity. The reason is that local and global weak continuities coincide in the conventional case.

**Proposition 5.** If $f$: $X \to Y$ is a locally weakly (Q, R)-continuous projection, then $f$ is weakly (Q, R)-continuous.

Indeed, let $f$ be a locally weakly (Q, R)-continuous projection and $V$ be an R-open set. Then $V$ is an R-neighborhood of some point $y$ from $Y$. As $f$ is a projection, there is a point $x$ from $X$ such that $f(x) = y$. By Definitions 10 a) and b), there is a neighborhood $Ox$ from $Q(x)$ such that $f(Ox)$ is a subset of $V$. This neighborhood $Ox$ is a Q-open set. As $V$ is an arbitrary R-open set, this concludes the proof of Proposition 5, demonstrating that the mapping $f$ is weakly (Q, R)-continuous.

Thus, we come to the following problems.

**Problem 1.** Is this result is true for arbitrary mapping, i.e., is a locally weakly (Q, R)-continuous mapping always weakly (Q, R)-continuous?

**Problem 2.** Is a weakly (Q, R)-continuous mapping always locally weakly (Q, R)-continuous?

In some cases, local and global weak (Q, R)-continuity do coincide.

**Proposition 6.** A projection $f: X \to Y$ is R-continuous if and only if it is locally R-continuous.

<u>Proof.</u> <u>Necessity</u>. Let us take an R-continuous projection $f: X \to Y$, a point $x$ from $X$, $y = f(x)$, and an R-neighborhood $Oy$ of $y$. As $f$ is R-continuous, $f^{-1}(Oy)$ is an open set. As the point $x \in f^{-1}(Oy)$, the set $f^{-1}(Oy)$ is a neighborhood of $x$ and $f(f^{-1}(Oy)) \subseteq Oy$. Thus, $f$ is locally R-continuous.

<u>Sufficiency</u> follows from Proposition 2.

Existence of two kinds of scaled continuity: weak scaled continuity and scaled continuity brings us to the following problems formulated in different terms in (Burgin, 2004).

**Problem 3.** Is weak (Q, R)-continuity weaker than (Q, R)-continuity?

**Problem 4.** Is weak (Q, R)-continuity at a point weaker than (Q, R)-continuity at a point?

The following example gives a positive solution to these problems.

**Example 15.** Let us take a two-dimensional Euclidean space $R^2$ and consider two spaces: the space $X$ that consists of all points $(x, y)$ such that $x$ belongs to the interval [0, 1] and $y$ belongs to the set {1, 2}, while the set $Y$ is the unit interval [0, 1] in the $x$-axis. The structure $Q_X(x)$ consists of $X$ and all connected open sets in $X$ that contain $x$. The structure $R_Y(y)$ consists of $Y$ and all connected open sets in $Y$ that contain $y$. The mapping $f: X \to Y$ is a natural projection of $X$ onto $Y$ defined by the formula $f(x, y) = (x, 0)$. Thus, the scale in $Y$ coincides with the standard topology of the interval [0, 1]. The mapping $f$ is weakly (Q, R)-continuous and weakly (Q, R)-continuous at any point of $X$, but it is neither (Q, R)-continuous at any point of $X$ nor (Q, R)-continuous globally as for any set $Z$ open in $Y$ its inverse image $f^{-1}(Z)$ is not connected and thus, does not belong to Q.

This example shows that the condition that Q is a C-structure is essential for coincidence of (Q, R)-continuity and weak (Q, R)-continuity as the system Q from this example is not a C-structure.

Structures Q and R play complimentary roles for fuzzy continuous mappings. The structure/scale R extends the scope of admissible functions, structure/scale Q restricts the scope of admissible functions. We can see (cf., for example, (Burgin, 2004)) that, in general there are much more R-continuous functions that continuous. However, the class of $(Q, T_X)$-continuous functions is much less than continuous ones.

Let $Y$ be a $T_1$, in particular, Hausdorff, topological space (Kelly, 1957) and Q be a P-structure, i.e., for any point $x$ from $X$, $Q_X(x)$ is the set of all open sets in $X$ that contain some chosen neighborhood $Ox$ of the point $x$. However, the class of $(Q, T_X)$-continuous functions is much less than .

**Theorem 3.** A mapping $f: X \to Y$ is weakly $(Q, T_X)$-continuous at a point $x$ from $X$ if and only if it is constant in some neighborhood of $x$.

<u>Proof</u>. Let us assume that $f$ is a weakly $(Q, T_X)$-continuous at the point $x$ mapping and there is a point $y$ in the chosen neighborhood $Ox$ such that $f(y) \neq f(x)$. As $Y$ is a $T_1$ topological space, there is a neighborhood $Of(x)$ of the point $f(x)$ such that $f(y)$ does not belong to $Of(x)$. Then by the definition of a weakly $(Q, T_X)$-continuous at the point $x$ mapping, there is a neighborhood $O_1x$ of the point $x$ for which $f(O_1x)$ is a subset of $Of(x)$ and $O_1x$ belongs to $Q_X(x)$. By the definition of the discontinuity structure Q, the set $O_1x$ contains $Ox$. Consequently, $f(y)$ belongs to $f(O_1x)$ and thus, belongs to $Of(x)$. This contradicts our assumptions and by principle of excluded middle shows that $f$ is constant in neighborhood $Ox$ of $x$.

Theorem is proved.

**Corollary 10.** A mapping $f: X \to Y$ is weakly $(Q, T_X)$-continuous if and only if it is constant in all connected components of $X$.

<u>Proof</u>. Let us assume that $Y$ is a connected component of $X$ and there are two points $x$ and $y$ in $Y$ such that $f(y) \neq f(x)$. By Proposition 3, the set $f^{-1}(f(x))$ is closed in X. At the same time, by Theorem 3, $f^{-1}(f(x))$ contains an open neighborhood for each of its points as $f$ is constant in some neighborhood of *x*. Consequently (Kelly, 1957), $f^{-1}(f(x))$ is an open subset of the set *Y*. This open and closed set $f^{-1}(f(x))$ does not coincide with *Y* because $f(y) \neq f(x)$. Consequently (Kelly, 1957), *Y* is not a connected subset of *X*. This contradicts our assumptions and by principle of excluded middle concludes the proof.

As Q is a C-structure, any weakly (Q, $T_X$)-continuous mapping is (Q, $T_X$)-continuous. It gives us the following result.

**Corollary 11.** A mapping *f*: *X*→ *Y* is (Q, $T_X$)-continuous at a point *x* from *X* if and only if it is constant in some neighborhood of *x*.

Let *Y* be a $T_1$ topological space and *X* be a metric space.

**Corollary 12.** For any number $a \in \mathbf{R}^{++}$, a mapping *f*: *X*→ *Y* is weakly ($Q_a$, $T_X$)-continuous at a point *x* from *X* if and only if it is constant in some neighborhood of *x*.

**Corollary 13.** For any number $a \in \mathbf{R}^{++}$, a mapping *f*: *X*→ *Y* is weakly ($Q_a$, $T_X$)-continuous if and only if it is constant in all connected components of *X*.

Let us consider three scaled topological spaces (*X*, $T_X$, Q), (*Z*, $T_Z$, P), and (*Y*, $T_Y$, R) and two mappings *f*: *X*→ *Y* and *g*: *Y*→ *Z*.

**Theorem 4.** If the mapping *f* is weakly (Q, R)-continuous at a point *x* (weakly (Q, R)-continuous) and the mapping *g* is weakly (R, P)-continuous at the point *f(x)* (weakly (R, P)-continuous), then the mapping *gf* is weakly (Q, P)-continuous at a point *x* (weakly (Q, P)-continuous, correspondingly).

Theorem 4 allows one to define a category **WSCTOP** of scaled topological spaces and their weakly fuzzy continuous mappings. Example 9 shows that categories **SCTOP** and **WSCTOP** do not coincide.

It is important to know relations between fuzzy and weakly fuzzy continuous mappings of the same topological spaces with different scales.

Let us consider two scaled topological spaces (*X*, $T_X$, $Q_X$) and (*X*, $T_X$, $P_X$).

**Definition 11** (Burgin, 2004). The scale $P_X$ is called *finer* than the scale $Q_X$ if for all *x* any element from $Q_X(x)$ contains as a subset some element from $P_X(x)$. The relation 'to be finer' is denoted by $Q_X \leq P_X$.

**Example 16.** The trivial scale $T_X$ (cf. Example 7) is finer than any other scale in $X$.

Let us consider how the relation to be finer influences relations between continuous mappings.

**Proposition 7.** a) For any scaled topological spaces $(X, T_X, Q)$, $(X, T_X, P)$, and $(Y, T_Y, R)$ if P is a C-structure that is finer than Q, then any (Q, R)-continuous (at $x$) mapping is (P, R)-continuous (at $x$).

b) For any scaled topological spaces $(X, T_X, Q)$, $(Y, T_Y, V)$, and $(Y, T_Y, R)$ if R is finer than V and either Q is a C-structure or R is a C-structure, then any (Q, R)-continuous (at $x$) mapping is (Q, V)-continuous (at $x$).

<u>Proof.</u> a. Let P be a C-structure that is finer than Q, $x$ is a point from $X$, and $f: X \to Y$ be a (Q, R)-continuous (at $x$) mapping. Then by Definition 9, for $y = f(x)$ and any neighborhood $Oy$ from $R(y)$, the set $f^{-1}(Oy)$ is an element of $Q(x)$. By Definition 11, $f^{-1}(Oy)$ contains some neighborhood $Ox$ from $P(x)$. As $P(x)$ is a filter, $f^{-1}(Oy)$ is also an element from $P(x)$. Consequently, the mapping $f$ is (P, R)-continuous (at $x$).

b. Let R be a C-structure that is finer than V, $x$ is a point from $X$, and $f: X \to Y$ be a (Q, R)-continuous (at $x$) mapping. Let us take some neighborhood $Oy$ from $V(y)$. The set $Oy$ contains a neighborhood $O_1 y$ from $R(x)$ because R is finer than V. At first, we consider the case when Q is a C-structure. By Definition 9, the set $f^{-1}(O_1 y)$ is an element of $Q(x)$. The set $f^{-1}(O_1 y)$ is a subset of $f^{-1}(Oy)$. As $Q(x)$ is a filter, $f^{-1}(Oy)$ is also an element from the neighborhood $Q(x)$. Consequently, the mapping $f$ is (Q, V)-continuous (at $x$).

When R is a C-structure, the set $Oy$ also belongs to $R(x)$ by properties of filters. Then by Definition 9, the set $f^{-1}(Oy)$ is also an element of $Q(x)$. Consequently, as $Oy$ is an arbitrary neighborhood of $y$, the mapping $f$ is (Q, V)-continuous (at $x$).

Proposition 7 is proved.

The next result is proved in a similar way.

**Proposition 8.** a) For any scaled topological spaces $(X, T_X, Q)$, $(X, T_X, P)$, and $(Y, T_Y, R)$ from $TQ \subseteq TP$ and $Q \subseteq P$ (as relations), it follows that any (Q, R)-continuous (at $x$) mapping is (P, R)-continuous (at $x$).

b) For any scaled topological spaces ($X$, $T_X$, $Q$), ($Y$, $T_Y$, $V$), and ($Y$, $T_Y$, $R$) from $TV \subseteq TR$ and $V \subseteq R$ (as relations), it follows that any ($Q$, $R$)-continuous (at $x$) mapping is ($Q$, $V$)-continuous (at $x$).

**Corollary 14.** For any scaled topological spaces ($Y$, $T_Y$, $V$) and ($Y$, $T_Y$, $R$) from $TV \subseteq TR$ and $V \subseteq R$ (as relations), it follows that any $R$-continuous (at $x$) mapping is $V$-continuous (at $x$).

**Corollary 15.** Any ($Q$, $R$)-continuous (at $x$) mapping $f: X \to Y$ is $R$-continuous (at $x$).

**Corollary 16.** Any continuous (at $x$) mapping $f: X \to Y$ is $R$-continuous (at $x$).

These results demonstrate that the concepts of $R$-continuity and ($Q$,$R$)-continuity are natural extensions of the concept of the conventional continuity.

In some cases, ($Q$, $R$)-continuity coincide with the conventional continuity. For example, let us take as the set $TR$ some base (Kuratowski, 1968) of the topology $T_Y$ and such $R$ that maps each point $x$ onto the set of all elements from $TR$ that contain the point $x$.

**Corollary 17.** A mapping $f: X \to Y$ is $R$-continuous (at $x$) if and only if it is continuous (at $x$).

Let ($X$, $T_X$, $Q$), ($Y$, $T_Y$, $R$), and ($Z$, $T_Z$, $P$) be scaled topological spaces.

**Proposition 9.** If a mapping $f: X \to Y$ is ($Q$, $R$)-continuous (at a point $x$ from $X$) and a mapping $g: Y \to Z$ is ($R$, $P$)-continuous (at a point $f(x)$), then the mapping $gf: X \to Z$ is ($Q$, $P$)-continuous (at a point $x$).

**Remark 5.** It is possible that a mapping $f: X \to Y$ is $R$-continuous (at a point $x$) and a mapping $g: Y \to Z$ is $P$-continuous (at a point $f(x)$), but the mapping $gf: X \to Z$ is not $P$-continuous (at a point $x$). Even if we take two $Q$-continuous (at a point $x$) mappings $f, g: X \to X$, their composition $gf$ might be not $Q$-continuous (at a point $x$). It is demonstrated by the following example.

**Example 17.** Let $X = [0,2]$; $a=1/10$; and $Q$ is equal to $Q_a$ that is restricted to this interval $[0,2]$. That is, $Q$ corresponds: to each $x$ from the interval $(1/10, 19/10)$ all intervals ($x + k$, $x - k$) with $1/10 < k$ and $0 \leq x - k < x + k \leq 2$; to each $x$ from the interval $[0, 1/10]$ all intervals [0, $x + k$) with 1/10 < $k$ and $x + k \leq 2$; to each $x$ from the interval [19/10, 2] all intervals ($x - k$, 2] with 1/10 < $k$ and $0 \leq x - k$.

We define $f: X \to X$ as follows. If $x$ belongs to the interval [0, 1), then $f(x) = 10x/11$, and $f(x) = x$ for all other $x$. In this case, $f$ has only one gap and this gap is less than 1/10. Consequently, $f$ is Q-continuous. However, $f^2$ has a gap that is equal to 21/121 because $f^2(1) = 100/121$ while for any $a > 1$, $f^2(a) = a$. As 21/121 > 1/10, the function $f^2$ is not Q-continuous.

A scale in the same topological space can be locally finer than another scale in the same space.

**Definition 12.** The scale $P_X$ is called *finer at a point x* than the scale $Q_X$ if any element from $Q_X(x)$ contains as a subset some element from $P_X(x)$. The relation 'to be finer at a point $x$' is denoted by $Q_X(x) \leq P_X(x)$.

It is possible that a scale in a topological space is finer at some point than another scale in the same space but the second scale is finer at another point than the first scale. The relation "to be globally finer" implies the relation "to be locally finer".

Let us consider two scaled topological spaces $(X, T_X, Q)$ and $(Y, T_Y, R)$, and a system $\mathbf{L} = \{R_i ; i \in I\}$ of scales in $Y$. We assume that $\mathbf{L}$ is a *base* for R at any point $y$ from $Y$, i.e., for any point $y$ from $Y$ and any neighborhood $U$ of $y$ from $R(y)$, there is a scale $R_i$ from $\mathbf{L}$ and a neighborhood $V$ of $y$ from $R_i(y)$, such that $V \subseteq U$, and the scale R is *finer* than any $R_i$ at any point $y$ from $Y$, i.e., for any point $y$ from $Y$, for any scale $R_i$ from $\mathbf{L}$, and any neighborhood $V$ of $y$ from $R_i(y)$, there is neighborhood $U$ of $y$ from $R(y)$ such that $U \subseteq V$.

**Theorem 5.** A mapping $f: X \to Y$ is weakly (Q, R)-continuous if the mapping $f$ is weakly $(Q, R_i)$-continuous for all $i \in I$.

<u>Proof.</u> <u>Necessity.</u> Let $f: X \to Y$ be a weakly (Q, R)-continuous mapping, $x$ be a point in the scaled topological space $(X, T_X, Q)$, and $V$ be a neighborhood of $f(x)$ from $R_i(f(x))$. As R is finer than $R_i$, there is neighborhood $U$ of $y$ from $R(y)$ such that $U \subseteq V$. By the definition of a weakly (Q, R)-continuous mapping, there is such neighborhood $Ox$ from $Q(x)$ that $f(Ox)$ is a subset of $U$. Consequently, $f(Ox)$ is a subset of $V$. As the neighborhood $V$ of $f(x)$ was arbitrary in $R_i$, this means that $f: X \to Y$ is a weakly (Q, $R_i$)-continuous mapping. As $R_i$ was arbitrary, this concludes the proof of necessity.

Sufficiency. Let $f: X \to Y$ be a weakly (Q, $R_i$)-continuous mapping for all $i \in I$, $x$ be a point in the scaled topological space (X, $T_X$, Q), and U be a neighborhood of $f(x)$ from R($f(x)$). As **L** is a base for R at any point $y$ from Y, there is a scale $R_i$ from **L** and a neighborhood V of $y$ from $R_i(y)$, such that $V \subseteq U$. By the definition of a weakly (Q, $R_i$)-continuous mapping, there is such neighborhood Ox from Q($x$) that $f(Ox)$ is a subset of V. Consequently, $f(Ox)$ is a subset of U. As the neighborhood U of $f(x)$ was arbitrary in $R_i$, this means that $f: X \to Y$ is a weakly (Q, R)-continuous mapping. This concludes the proof of necessity, and thus, the proof of the theorem.

In a similar way, we can prove the following local form of Theorem 5.

Let us consider a point $x$ in the scaled topological space (X, $T_X$, Q). Assume that **L** is a *base* for R at the point $f(x)$, i.e., for any neighborhood U of $f(x)$ from R($f(x)$), there is a scale $R_i$ from **L** and a neighborhood V of $f(x)$ from $R_i(f(x))$ such that $V \subseteq U$, and the scale R is *finer* at the point $f(x)$, i.e., for any scale $R_i$ from **L** and a neighborhood V of $f(x)$ from $R_i(x)$ there is a neighborhood U of $x$ from R($f(x)$) such that $U \subseteq V$.

**Theorem 6.** A mapping $f: X \to Y$ is weakly (Q, R)-continuous at a point $x$ (weakly (Q, R)-continuous) if the mapping $f$ is weakly (Q, $R_i$)-continuous (weakly (Q, $R_i$)-continuous, correspondingly) at the point $x$ for all $i \in I$.

### 4. Conclusion

Thus, it is demonstrated how a basic topological concept such as continuity is extended to scaled or fuzzy continuity in scaled topological spaces. The goal of this extension is to adjust topology to demands of contemporary multiscale physics and numerical simulation of processes going at different scales. This study is done in discontinuous or scaled topology, which is a new field, in which methods and constructions of the classical topology are scaled.

This combination of scales and topological methods gives birth to the name "scaled topology."

Results obtained in this paper make possible to obtain classical results for continuous mappings as direct corollaries, which do not demand additional proofs. This demonstrates that scaled topology is a natural extension of the classical topology, which is aimed at a more adequate representation of real-life situations such as computation or measurement of topological characteristics.

Scaled topology in real linear spaces has applications to computational mathematics and numerical methods (Burgin and Westman, 2000). Consequently, it is interesting to study connections between scaled topology and interval analysis, which stemmed from problems of computations (Moore, 1966).

Scalable topology allows one to introduce scalability to digital topology and to develop further this field bringing more flexibility to its methods (cf., Rosenfeld, 1979; Smyth, 1995). According to Rosenfeld, digital topology is the use of known facts in topology or graph theory to understand the topological properties of the finite images that come up in computer science, well enough to write certain algorithms for them.

There are related problems of continuity in software engineering, particularly in testing and certifying programs (Hamlet, 2002). Scalable topology allows one to define rigorously software continuity that may be used for different applications and continuity of software specifications.

Now topology is developed not only in sets but also in the context of abstract categories and functors (cf., for example, (Johnstone, 1977)). Thus, it might be interesting to build scalable topology in a more abstract, categorical setting.

At the same time, it is necessary to remark that physics utilizes mostly the topology of manifolds, which is more concrete that general topological spaces considered in this work. However, the scalability technique developed here allows one to build scales in the topology of manifolds.

One more direction for scalable topology applications is the further development of discrete topology. Scalability allows one to study discrete spaces as subspaces of classical continuous spaces and utilize several scales in the same space.